\documentclass[12pt]{article}
\usepackage{amsmath}
\usepackage{amsfonts}
\usepackage{amssymb}
\usepackage{amsthm}
\usepackage{hyperref}

\usepackage{graphicx}

\def\bee{\begin{equation}}
\def\eee{\end{equation}}

\begin{document}

\thispagestyle{empty}
\centerline{}
\bigskip
\bigskip
\bigskip
\bigskip
\bigskip
\bigskip
\centerline{\Large\bf Are the Stieltjes  constants irrational? }

\bigskip

\centerline{\Large\bf  Some computer experiments}
\bigskip
\bigskip
\bigskip

\begin{center}
{\large \sl Krzysztof D. Ma{\'s}lanka $^1$,  Marek Wolf $^2$}\\*[5mm]
\begin{center}
$^1$Institute for the History of Science of the Polish Academy of Sciences,\\
ul. Nowy {\'S}wiat 72, pok. 9,  00-330 Warsaw, e-mail:  krzysiek2357@gmail.com
\end{center}

\begin{center}
$^2$Cardinal  Stefan  Wyszynski  University, Faculty  of Mathematics and Natural Sciences. College of Sciences,\\
ul. W{\'o}ycickiego 1/3,   PL-01-938 Warsaw,   Poland, e-mail:  m.wolf@uksw.edu.pl
\end{center}

\bigskip
\end{center}

\bigskip
\bigskip

\begin{center}
{\bf Abstract}\\
\bigskip
\begin{minipage}{12.8cm}
Khnichin's theorem is a surprising and still relatively little known result.
It can be used as a specific criterion for determining whether or not any
given number is irrational. In this paper we apply this theorem    as well as the Gauss--Kuzmin theorem
to several   thousand high precision (up to  more  than  53000 significant digits) initial Stieltjes
constants $\gamma _{n}$, $n=0,1,...,5000$ in order to confirm that, as is
commonly believed, they are irrational numbers  (and even transcendental).  We  study also the normality of these important  constants.
\end{minipage}
\end{center}

\bigskip\bigskip

\bibliographystyle{abbrv}

\section{Introduction}

The famous zeta function $\zeta (s)$ discovered by L. Euler in 1737   and published in 1744    \cite{Euler_1737}
as a function of  real variable was investigated by G. F. B. Riemann  in the complex
domain in his famous memoir submitted in 1859 to the Prussian Academy
\cite{Riemann-1859}. It is defined as:%
\begin{equation}
\zeta (s)=\sum\limits_{n=1}^{\infty }\frac{1}{n^{s}}\qquad \Re(s)>1.
\label{Zeta}
\end{equation}%
It is divergent in the most interesting area of the complex plane, i.e., in
the so called critical strip $0\leq \Re(s)\leq 1$ where all complex
zeros of zeta lie. However, as was shown by Riemann, the definition
\eqref{Zeta} does contain information about the zeta function on the entire
complex plane but the process of analytic continuation must be used in order
to reveal global behavior of this function.   In fact Riemann in his paper
analytically continued \eqref{Zeta}  to the whole
complex plane except   $s=1$ by means of the following contour integral:
\bee
\zeta(s)=\frac{\Gamma(-s)}{2\pi i}\underset{\mathcal{P}}\int \frac{(-x)^s}{e^x-1} \frac{dx}{x},
\eee
\noindent where the integration is performed along the  following path  ${\mathcal{P}}$:

\vskip 0.3cm
\begin{picture} (70,120)(0,0) \thicklines
  \put(180,40){\vector(0,1){80}} \put(120,80){\vector(1,0){130}}
\put(180,80){\oval(40,40)[l]}
\put(180,85){\oval(30,30)[rt]}
\put(180,75){\oval(30,30)[rb]}

\put(195,75){\vector( 1,0){45}}
\put(240,85){\vector(-1,0){45}}
\end{picture}.

\vskip -1.2cm
\noindent   Till now dozens of  integrals   and series representing the  $\zeta(s)$
function are known,  for collection of  such formulas see for example  the entry
{\it Riemann Zeta Function} in \cite{Weisstein} and references cited therein  and  \cite{Milgram-2013}.

Another representation of this function is given by a power series where
appear certain constants $\gamma _{n}$. These constants are essentially
coefficients of the Laurent series expansion of the zeta function around its
only simple pole at $s=1$:%
\begin{equation}
\zeta (s)=\frac{1}{s-1}+\sum\limits_{n=0}^{\infty }\frac{\left( -1\right)^{n}}{n!}\gamma _{n}\left( s-1\right) ^{n}
\label{ZetaExpansion}
\end{equation}%
Primary definition of these fundamental constants was found by Th. J.
Stieltjes and presented in a letter to Ch. Hermite dated June 23, 1885
\cite[letters no. 71--74]{Hermite_Charles_(1822-1901)_Correspondance_1905}

\begin{equation}
\gamma _{n}=\underset{m\rightarrow \infty }{\lim }\left[ \left(
\sum\limits_{k=1}^{m}\frac{ (\ln k)^{n}}{k}\right) -\frac{( \ln m)^{n+1}}{n+1}\right].
\label{gamma}
\end{equation}%
(When $n=0$ the     numerator of the fraction     in the first summand in (\ref{gamma})
is formally $0^{0}$ which is  taken to be $1$.)

Effective numerical computing of the constants
$\gamma_n$ is quite a challenge because the formulas \eqref{gamma}
converge extremely slowly.  Even when $n=0$,
which corresponds to the well-known Euler-Mascheroni constant $\gamma _{0}$,
in order to obtain just $10$ accurate digits one has to sum up exactly $%
12366 $ terms whereas in order to obtain $10000$ digits (which is indeed
required in some applications) one would have to sum up unrealistically
large number of terms: nearly $5\cdot 10^{4342}$ which is of course far
beyond capabilities of the present day computers.
However, various fast algorithms were found to efficiently compute specific value of the
zeroth Stieltjes constant $\gamma_0$, i.e. the fundamental Mascheroni-Euler constant,
see e.g.  \cite{Sweeney_1963}, \cite{Brent_McMillan_1980}.     For $n>0$ the situation is
still worse. Therefore we have to seek for other faster algorithms. In 1992
J. B. Keiper \cite{Keiper_1992} published an effective algorithm based on numerical
quadrature of certain integral representation of the zeta function and
alternating series summation using Bernoulli numbers. Keiper's algorithm was
later implemented in widely used program \textit{Mathematica}.\ An efficient
but rather complicated method based on Newton-Cotes quadrature has been
proposed by R. Kreminski in 2003 \cite{Kreminski_2002}. Quite recently F.  Johansson
presented particularly efficient method \cite{Johansson_2014}.

In the Appendix at the end of the present paper yet another method of
computing Stieltjes constants will be described which is perhaps
not as efficient as Johansson's approach, yet it is by far more simple and
it may be easily and quickly used in practical calculations for obtaining
$\gamma_{n}$ up to $n\sim 10000$ with accuracy $\sim 50000$ significant
digits.

We proceed as follows. First, we use the algorithm presented in the Appendix to calculate 5000 $\gamma_n$  with accuracies ranging from about
53000 significant digits ($\gamma_0$) to about 24000 digits ($\gamma_{5000}$). Having these numbers we intend to provide an argument in
favor of their irrationality. Then we consider the question of their normality, as real expansions in the base equal 10.
Finally, in Sect. 3, we develop $\gamma_n$'s into continuous fractions and next use the remarkable theorems due to Khinchin, L{\'e}vy   and Gauss--Kuzmin.
Obtained results support the common opinion that   $\gamma_n$ are indeed irrational.

\section{Normality}

Let us recall that  a number $r$ is normal  in base $b$  if  each finite string of $k$ consecutive digits appears in
this expansion with  asymptotic frequency  $b^{-k}$.  In the usual decimal
base we have  that  each digit $0,1,2, \ldots, 9$  appears in  the  expansion of  the number $r$ with limiting   frequency 0.1,
each 2--digits  string  $00,  01, \ldots, 99$  appears with density 0.01.   Having the first 5000 Stieltjes  constants
with accuracies as  described earlier we checked that each digit $0,1, 2, \ldots, 9$  appears almost exactly with frequency 0.1.
It is difficult to represent this  $5000\times  10$  data points  in one plot.  In the Fig. 1 we employed the following artifice:  the frequency
$h_n(0)$ of   appearance of  digit 0 in the Stieltjes  constant $\gamma_n$  is plotted at $x$--axis value $n$  with the $y$ value $0.1-h_ n(0)$,
i.e. the distance from the  expected value 0.1,  which in this case of $a=0$  should be  around 0.1.    In general,  the  frequency  $h_n(a)$
of appearance of digit $a$ in the Stieltjes  constant
$\gamma_n$  is plotted   with the $y$ value $a\times  0.1+(0.1-h_n(a))$.  We calculated also density of
100 strings of two digits $00, 01, \ldots, 99$  for all 5000  Stieltjes  constants  $\gamma_n$. Now the result consisted of half a million
points,  what is impossible to represent on the plot. Instead,  in the Table I we present for each pattern of digits $ab$ the maximal
difference between  calculated  frequency of appearance and the expected value of 0.01 and the  number $n$ of the Stieltjes  constant  $\gamma_n$
for which this  discrepancy appeared.  The difference between the actual computed value
of the frequency of two digits patterns and the
expected value 0.01 was  typically of a few  percents.

\begin{figure}[h]
\begin{center}
\includegraphics[width=0.8\textwidth, angle=0]{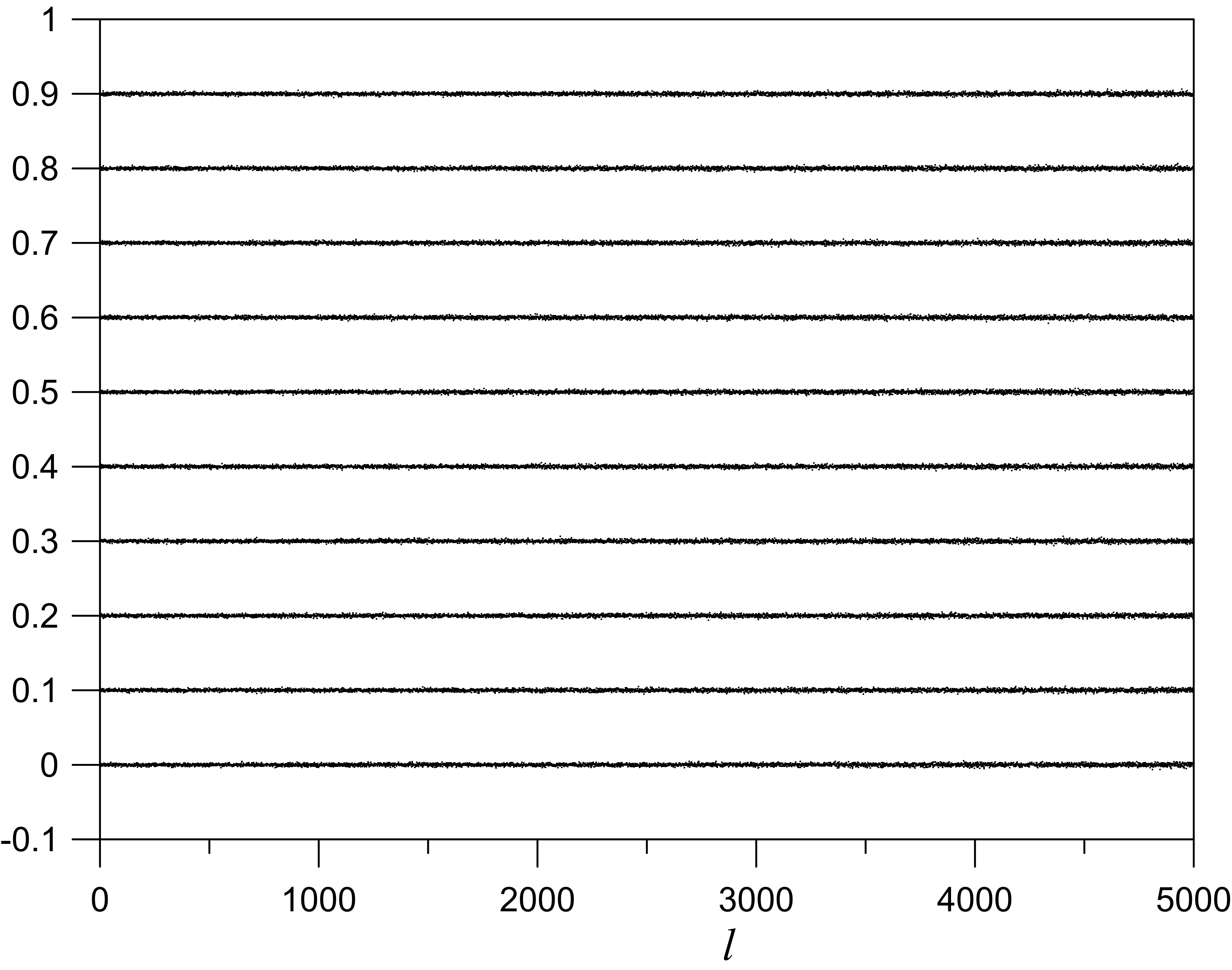} \\
\vspace{-1.5cm}
\vspace{1.7cm}Fig.1  The plot of the differences between 0.1 and  actual frequencies of digits $0, 1, \dots, 9$  for all 5000
 Stieltjes  constants.  The data for  digit  $a$  is plotted  at $y$  value $a\times 0.1$  for clarity. \\
\end{center}
\end{figure}

\section{Continued fractions  expansions }

Continued fractions often reveal various profound and unexpected properties of irrational
numbers that are normally hidden in their traditional decimal (or other basis) notation,  see e.g. \cite{Martin_2004}.

In this Section  we are going to exploit three  facts about the continued fractions: the
existence of the   Khinchin constant, Khinchin--L{\'e}vy constant   and   the Gauss--Kuzmin  distribution,
see e.g. \cite[chapter III,  \S 15] {Khinchin},\cite[\S 1.8,  \S2.17]{Finch},  to support the irrationality of Stieltjes  constants  $\gamma_n$.
The paper \cite{Brent_1977}  presents  the regular continued fraction for the  Euler's--Mascheroni  constant $\gamma_0$.  Let
\bee
r = [a_0(r); a_1(r), a_2(r), a_3(r), \ldots]=
a_0(r)+\cfrac{1}{a_1(r) + \cfrac{1}{a_2(r) + \cfrac{1}{a_3(r)+\ddots}}}
\eee
be the continued fraction expansion of the real number $r$, where $a_0(r)$ is an integer and all denominators  $a_k(r)$ (``partial  quotients'')
with $k\geq 1$  are positive integers.  Let us remark  that  rational numbers  have finite  number of  coefficients  $a_k$.
Khinchin has proved \cite{Khinchin},  see also \cite{Ryll-Nardzewski1951}, that  limits of geometrical means of $a_k(r)$ are the same
for almost all real $r$:
\bee
\lim_{l\rightarrow \infty} \big(a_1(r) \ldots a_l(r)\big)^{\frac{1}{l}}=
\prod_{m=1}^\infty {\left\{ 1+\frac{1}{ m(m+2)}\right\}}^{\log_2 m} \equiv K_0=  2.685452001\dots   ~~.
\label{Khinchin}
\eee
The Lebesgue measure of (all) the exceptions is zero  and include {\it rational numbers}, quadratic irrationals and
some irrational numbers too, like for example the  Euler constant $e=2.7182818285\ldots$ for which the limit (\ref{Khinchin}) is
infinity.

\vskip -2cm
\begin{center}
{\centerline{ \bf Table 1}}
In the columns  A, C, E and G the two digits patterns are given, in the columns B, D, F  and H the maximal differences between 0.01 and
the frequency  that a given pattern $ab$, $a, b=0,1,\ldots, 9$  appears  among the digits of the   $\gamma_n$,  $n=1,2,\ldots, 5000$. \\
\end{center}
\begin{center}
{\small
\begin{tabular}{|c|c|c|c|c|c|c|c|} \hline
A & B & C & D & E & F & G & H\\ \hline
00  & $  2.4914\times 10^{-3} $ &  $ 25  $  &  $ 1.7898\times 10^{-3} $  &  $ 50 $ &  $ 2.0046\times 10^{-3}  $  &  $75  $  &  $ 1.9586\times 10^{-3} $\\ \hline
01  & $  2.0114\times 10^{-3} $ &  $ 26  $  &  $ 2.3187\times 10^{-3} $  &  $ 51 $ &  $ 2.1064\times 10^{-3}  $  &  $76  $  &  $ 2.0058\times 10^{-3}$  \\  \hline
02  & $  2.0771\times 10^{-3} $ &  $ 27  $  &  $ 2.0847\times 10^{-3} $  &  $ 52 $ &  $ 2.2251\times 10^{-3}  $  &  $77  $  &  $ 2.1520\times 10^{-3}$  \\  \hline
03  & $  2.3235\times 10^{-3} $ &  $ 28  $  &  $ 2.5891\times 10^{-3} $  &  $ 53 $ &  $ 2.2773\times 10^{-3}  $  &  $78  $  &  $ 2.1413\times 10^{-3}$  \\  \hline
04  & $  1.8466\times 10^{-3} $ &  $ 29  $  &  $ 2.1732\times 10^{-3} $  &  $ 54 $ &  $ 1.9028\times 10^{-3}  $  &  $79  $  &  $ 2.2307\times 10^{-3}$  \\  \hline
05  & $  1.9006\times 10^{-3} $ &  $ 30  $  &  $ 1.9310\times 10^{-3} $  &  $ 55 $ &  $ 2.2080\times 10^{-3}  $  &  $80  $  &  $ 1.8309\times 10^{-3}$  \\  \hline
06  & $  1.8525\times 10^{-3} $ &  $ 31  $  &  $ 2.0466\times 10^{-3} $  &  $ 56 $ &  $ 2.4565\times 10^{-3}  $  &  $81  $  &  $ 2.1083\times 10^{-3}$  \\  \hline
07  & $  2.4075\times 10^{-3} $ &  $ 32  $  &  $ 2.0625\times 10^{-3} $  &  $ 57 $ &  $ 1.8966\times 10^{-3}  $  &  $82  $  &  $ 1.8493\times 10^{-3}$  \\  \hline
08  & $  2.4080\times 10^{-3} $ &  $ 33  $  &  $ 2.1236\times 10^{-3} $  &  $ 58 $ &  $ 1.9259\times 10^{-3}  $  &  $83  $  &  $ 2.1614\times 10^{-3}$  \\  \hline
09  & $  2.0118\times 10^{-3} $ &  $ 34  $  &  $ 1.9970\times 10^{-3} $  &  $ 59 $ &  $ 2.0112\times 10^{-3}  $  &  $84  $  &  $ 2.3112\times 10^{-3}$  \\  \hline
10  & $  2.1949\times 10^{-3} $ &  $ 35  $  &  $ 2.2988\times 10^{-3} $  &  $ 60 $ &  $ 1.9846\times 10^{-3}  $  &  $85  $  &  $ 2.6315\times 10^{-3}$  \\  \hline
11  & $  2.3476\times 10^{-3} $ &  $ 36  $  &  $ 2.1588\times 10^{-3} $  &  $ 61 $ &  $ 1.9017\times 10^{-3}  $  &  $86  $  &  $ 1.9200\times 10^{-3}$  \\  \hline
12  & $  1.8161\times 10^{-3} $ &  $ 37  $  &  $ 2.2839\times 10^{-3} $  &  $ 62 $ &  $ 1.9813\times 10^{-3}  $  &  $87  $  &  $ 2.1604\times 10^{-3}$  \\  \hline
13  & $  1.9746\times 10^{-3} $ &  $ 38  $  &  $ 1.9860\times 10^{-3} $  &  $ 63 $ &  $ 2.3341\times 10^{-3}  $  &  $88  $  &  $ 2.4448\times 10^{-3}$  \\  \hline
14  & $  2.3346\times 10^{-3} $ &  $ 39  $  &  $ 2.1897\times 10^{-3} $  &  $ 64 $ &  $ 2.2752\times 10^{-3}  $  &  $89  $  &  $ 2.3153\times 10^{-3}$  \\  \hline
15  & $  2.1317\times 10^{-3} $ &  $ 40  $  &  $ 2.1021\times 10^{-3} $  &  $ 65 $ &  $ 1.9558\times 10^{-3}  $  &  $90  $  &  $ 1.8766\times 10^{-3}$  \\  \hline
16  & $  1.8801\times 10^{-3} $ &  $ 41  $  &  $ 2.2182\times 10^{-3} $  &  $ 66 $ &  $ 2.3915\times 10^{-3}  $  &  $91  $  &  $ 2.2997\times 10^{-3}$  \\  \hline
17  & $  1.8627\times 10^{-3} $ &  $ 42  $  &  $ 2.1976\times 10^{-3} $  &  $ 67 $ &  $ 2.3017\times 10^{-3}  $  &  $92  $  &  $ 2.1946\times 10^{-3}$  \\  \hline
18  & $  2.0085\times 10^{-3} $ &  $ 43  $  &  $ 1.9233\times 10^{-3} $  &  $ 68 $ &  $ 2.1579\times 10^{-3}  $  &  $93  $  &  $ 1.8714\times 10^{-3}$  \\  \hline
19  & $  2.3663\times 10^{-3} $ &  $ 44  $  &  $ 2.5452\times 10^{-3} $  &  $ 69 $ &  $ 1.8103\times 10^{-3}  $  &  $94  $  &  $ 1.8551\times 10^{-3}$  \\  \hline
20  & $  1.8711\times 10^{-3} $ &  $ 45  $  &  $ 1.9193\times 10^{-3} $  &  $ 70 $ &  $ 2.0240\times 10^{-3}  $  &  $95  $  &  $ 2.7646\times 10^{-3}$  \\  \hline
21  & $  2.0741\times 10^{-3} $ &  $ 46  $  &  $ 1.9071\times 10^{-3} $  &  $ 71 $ &  $ 1.9349\times 10^{-3}  $  &  $96  $  &  $ 1.9379\times 10^{-3}$  \\  \hline
22  & $  2.2366\times 10^{-3} $ &  $ 47  $  &  $ 2.1403\times 10^{-3} $  &  $ 72 $ &  $ 1.9635\times 10^{-3}  $  &  $97  $  &  $ 2.0152\times 10^{-3}$  \\  \hline
23  & $  2.2588\times 10^{-3} $ &  $ 48  $  &  $ 1.9612\times 10^{-3} $  &  $ 73 $ &  $ 1.9174\times 10^{-3}  $  &  $98  $  &  $ 1.9536\times 10^{-3}$  \\  \hline
24  & $  2.3669\times 10^{-3} $ &  $ 49  $  &  $ 1.9473\times 10^{-3} $  &  $ 74 $ &  $ 1.9815\times 10^{-3}  $  &  $99  $  &  $ 2.0863\times 10^{-3}$  \\  \hline
\end{tabular}
}
\end{center}

The constant $K_0$ is called the Khinchin constant, see e.g.  \cite[\S 1.8]{Finch}.
If the quantities
\bee
K(r; l)=\big(a_1(r) a_2(r) \ldots a_l(r)\big)^{\frac{1}{l}}
\label{Kny}
\eee
for a given number $r$ are close to $K_0$ we can regard it as an indication
that $r$ is irrational.

\begin{figure}[h]
\begin{center}
\includegraphics[width=0.8\textwidth, angle=0]{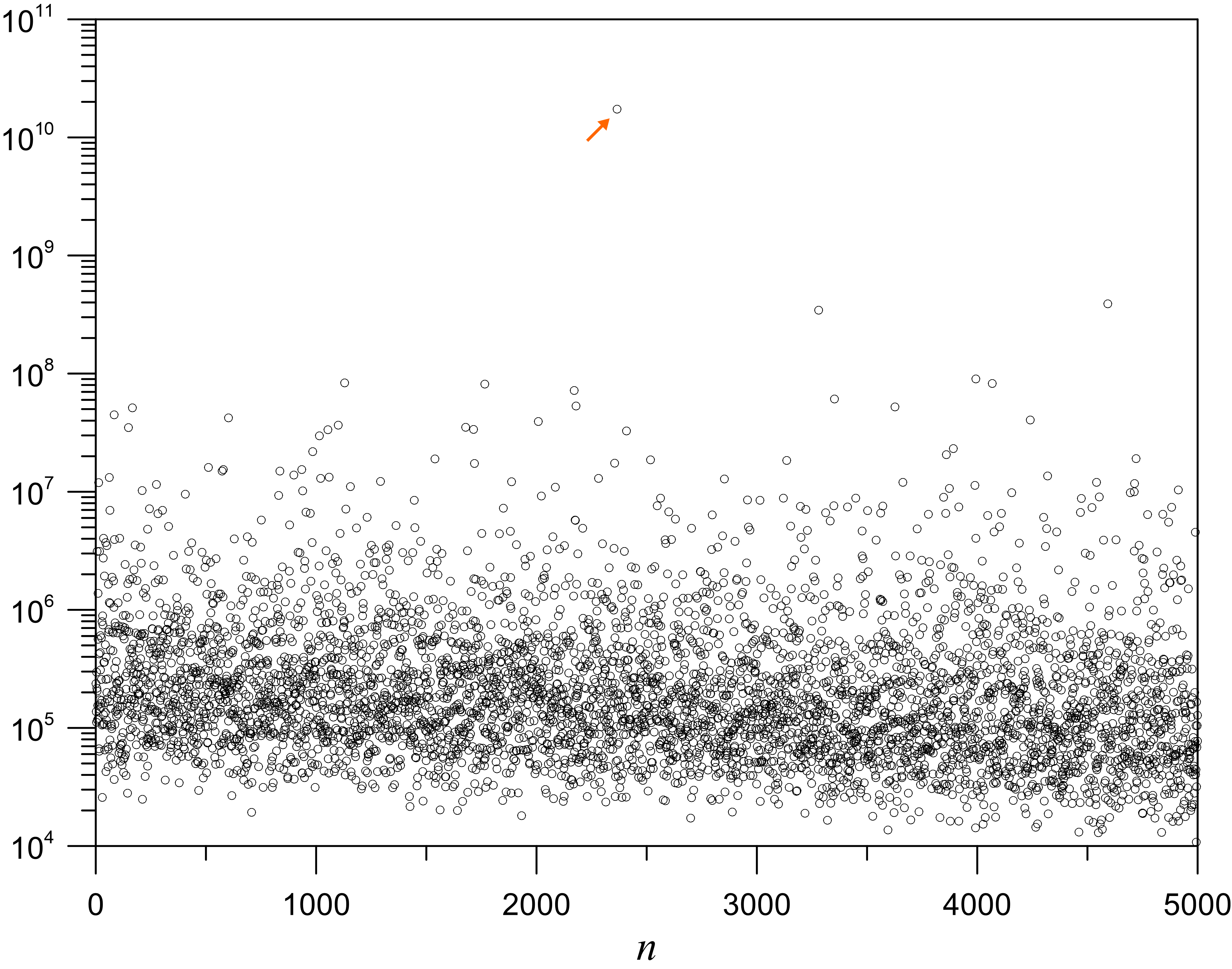} \\
\vspace{-1.5cm}
\vspace{1.7cm}Fig.2  The plot of maximal $a_k(n)$  for $n=1,2,3, \ldots, 5000$.     \\
\end{center}
\end{figure}

We developed the fractional parts of Stieltjes  constants (in Sect.2,   investigating the normality,  we used the whole  number,  e.g.
$\gamma_{62}=111670.9578149410793387893\ldots$  and we use in this section only digits after the decimal dot)  using  built in  PARI/GP   \cite{PARI2}
the  function  \verb"contfrac"$(r,\{nmax\})$ which creates the row vector ${\bf a}(r)$ whose
components are the denominators $a_k(r)$  of the continued fraction
expansion of $r$, i.e.  ${\bf a}=[a_0(r); a_1(r), \dots,a_l(r)]$ means that
\bee
r \approx
a_0(r)+\cfrac{1}{a_1(r) + \cfrac{1}{a_2(r) + \cfrac{1}{\ddots\cfrac{1}{a_l(r)}}}}
\eee
The parameter $nmax$ limits the number of terms $a_{nmax}(r)$; if it is omitted
the expansion stops with a declared precision
of representation of real number  $r$  at the last significant partial quotient:
the values of the convergents  ${P_k(r)/Q_k(r)}$
\bee
\frac{P_k(r)}{Q_k(r)}=a_0(r)+\cfrac{1}{a_1(r) + \cfrac{1}{a_2(r) + \cfrac{1}{a_3(r)+\ddots +\cfrac{1}{a_k}}}}
\eee

approximate the value of $r$  with accuracy
at least $1/Q_k^2$
\cite[Theorem 9, p.9]{Khinchin}:
\bee
\left|r - \frac{P_k(r)}{Q_k(r)} \right |<\frac{1}{Q_k^2(r)},
\label{error}
\eee
hence  when $1/{Q_k^2}$  is  smaller than  the accuracy of the number $r$  the process stops.

We checked that the PARI precision set to
{\tt $\backslash$p  120000} digits  is sufficient in the  sense that scripts with larger precision generated exactly
the same results: the rows  ${\bf a}(\gamma_n)$ obtained  with accuracy 140000 digits were the same for
all $n$  as those  obtained for accuracy 120000 and the continued fractions with  accuracy set to 100000 digits
had  different denominators $a_k(\gamma_n)$.   The number of partial  quotients $a_k$  varied  from  over 110000  for  initial   Stieltjes
constants  to  48027 for  $\gamma_{5000}$,   i.e.  the value of $l(n)$  was  roughly  2 times  the number of digits  in the expansion of $\gamma_n$.
However, there have been cases of extremely large values of  partial quotients.   The largest  was $a_{13034}=17399017050$  for  $\gamma_{2366}$,
marked by the red arrow at the top in Fig. 2.

\begin{figure}[h]
\begin{center}
\includegraphics[width=0.8\textwidth, angle=0]{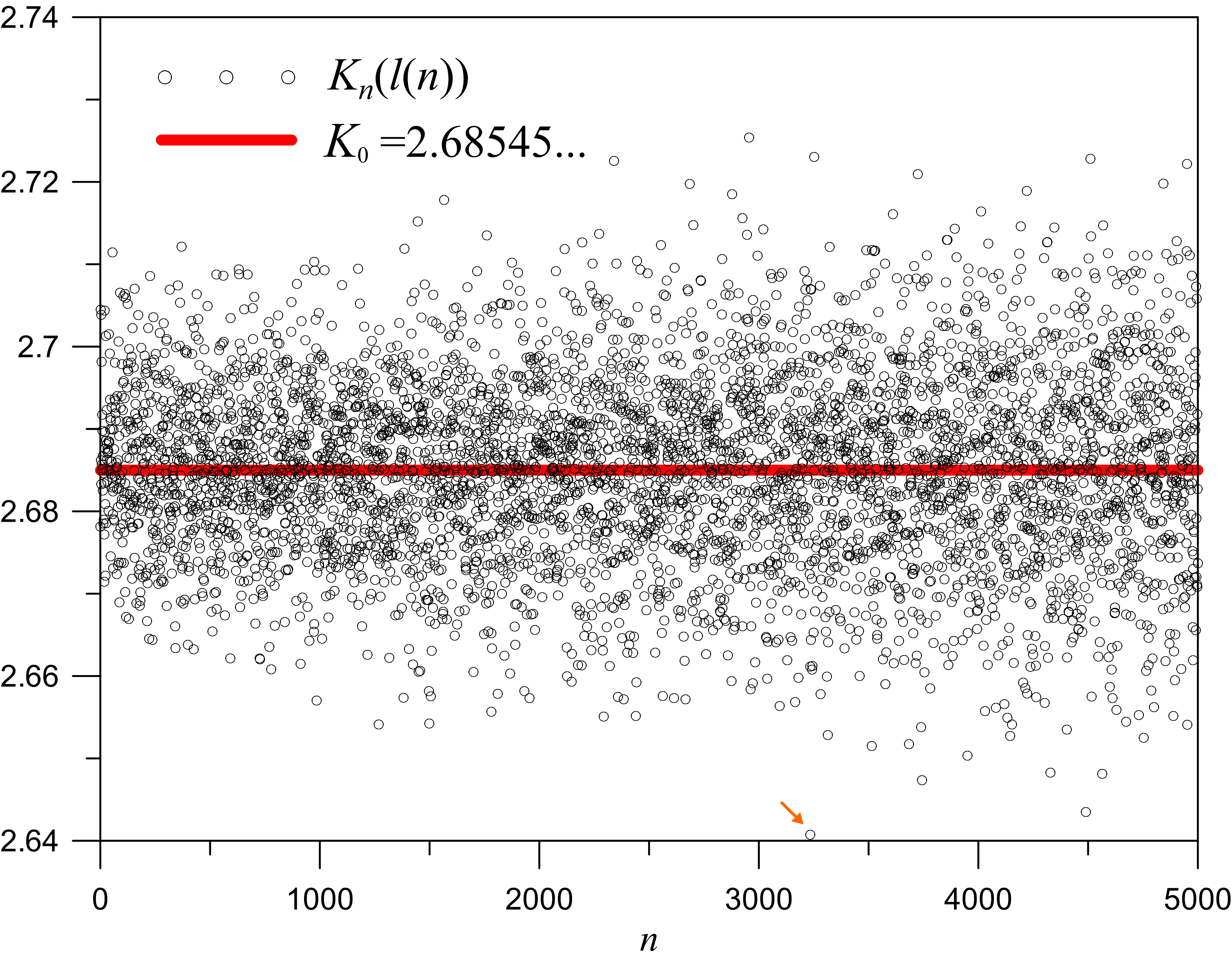} \\
\vspace{-1.5cm}
\vspace{1.7cm}Fig.3  The plot of $K_l(n(l))$ for $n=1,2,3, \ldots, 5000$.  There are 384  points  closer  to $K_0$ than $0.001$
and  30 points  closer  to $K_0$ than $0.0001$.
The largest  value of  $|K_0-K_n(l(n))|$ is $4.47\times 10^{-2}$ and it occurred for the Stieltjes constant  number $n=3235$ (marked with the red
arrow),  the smallest  value of  $|K_0-K_n(l(n))|$ is  $1.02\times 10^{-5}$    and it occurred for $\gamma_{1563}$.   \\
\end{center}
\end{figure}

With the precision set  to 120000 digits   we have expanded  each $\gamma_n$, $n=1, 2, \ldots 5000$    into  its  the continued fractions
($\doteq$   means ``approximately equal'')
\bee
\gamma_n \doteq [a_0(n); a_1(n), a_2(n), a_3(n), \ldots, a_{l(n)}(n)]\equiv {\bf a}(n)
\label{g-cfr}
\eee
without specifying the  parameter $nmax$, thus the length of the vector ${\bf a}(n)$  depended on $\gamma_n$
and it turns out that the number $l(n)$  of denominators  was contained between 53000 for  Stieltjes constants  with  index around 5000 and  110000
for gammas  with smallest  index $n$.  
The value of the product  $a_1 a_2 \ldots a_{l(n)} $  was typically of the
order $10^{47000}$  for  beginning  Stieltjes  constants to  $10^{23000}$  for the last  $\gamma_n$'s.  It means that,  if these  Stieltjes
constants are rational numbers $P/Q$   then $Q$  are larger then those  big numbers, for justification   see  e.g. \cite[Theorems  16, 17]{Khinchin}.
Next for each $n$ we have calculated the geometrical means:
\bee
K_n(l(n))= \left( \prod_{k=1}^{l(n)} a_k(n) \right)^{1/l(n)}.
\eee
The results are presented in the Fig.3. Values of $K_n(l(n))$ are scattered around the
red line representing $K_0$.  To gain some insight into the rate of convergence
of $K_n(l(n))$ we have plotted in the Fig. 4 the number of sign changes $S_K(n)$ of
$K_n(m)-K_0$ for each $n$ when $m=100, 101, \ldots l(n)$, i.e.
\bee
S_K(n)={\rm number ~ of ~ such ~{\it  m}~ that} ~~~(K_n(m+1)-K_0)(K_n(m)-K_0)<0.
\eee
The largest $S_K(n)$ was  961 and it occurred for the  $\gamma_{1175}$  and for 124
gammas  there were no sign changes at all. 
It is well known that the convergence to Khinchin's constant is very slow.   In the Fig.4   for each $\gamma_n$  we  present the   closest  to the
Khnichin constant $K_0$  value of  the ``running''  geometrical  means
\bee
K_n(m)=  \left( \prod_{k=1}^{m} a_k(n) \right)^{1/m},  ~~~~ m=100, 101, \ldots,  l(n).
\eee

\begin{figure}[h]
\begin{center}
\includegraphics[width=0.8\textwidth, angle=0]{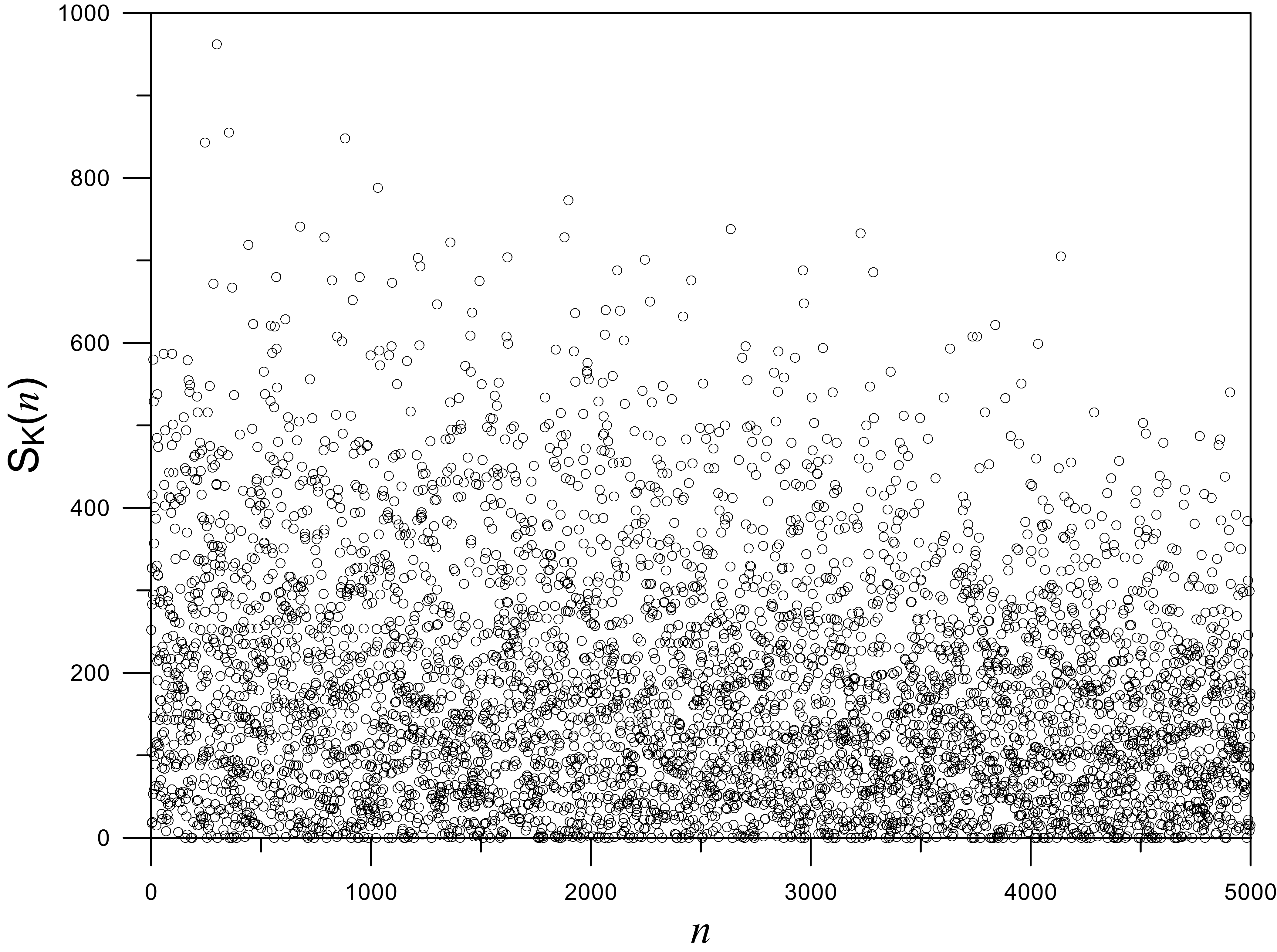} \\
\vspace{-1.5cm}
\vspace{1.7cm}Fig.4   The number of sign changes $S_K(n)$   for each  $n$,  i.e.  the number of  such $m$ that $(K_n(m+1)-K_0)(K_n(m)-K_0)<0$
(the initial transient values  of $m$ were skipped---  sign changes were detected for $m=100, 101, \ldots l(n)$).  \\
.   \\
\end{center}
\end{figure}

\begin{figure}[h]
\begin{center}
\includegraphics[width=0.8\textwidth, angle=0]{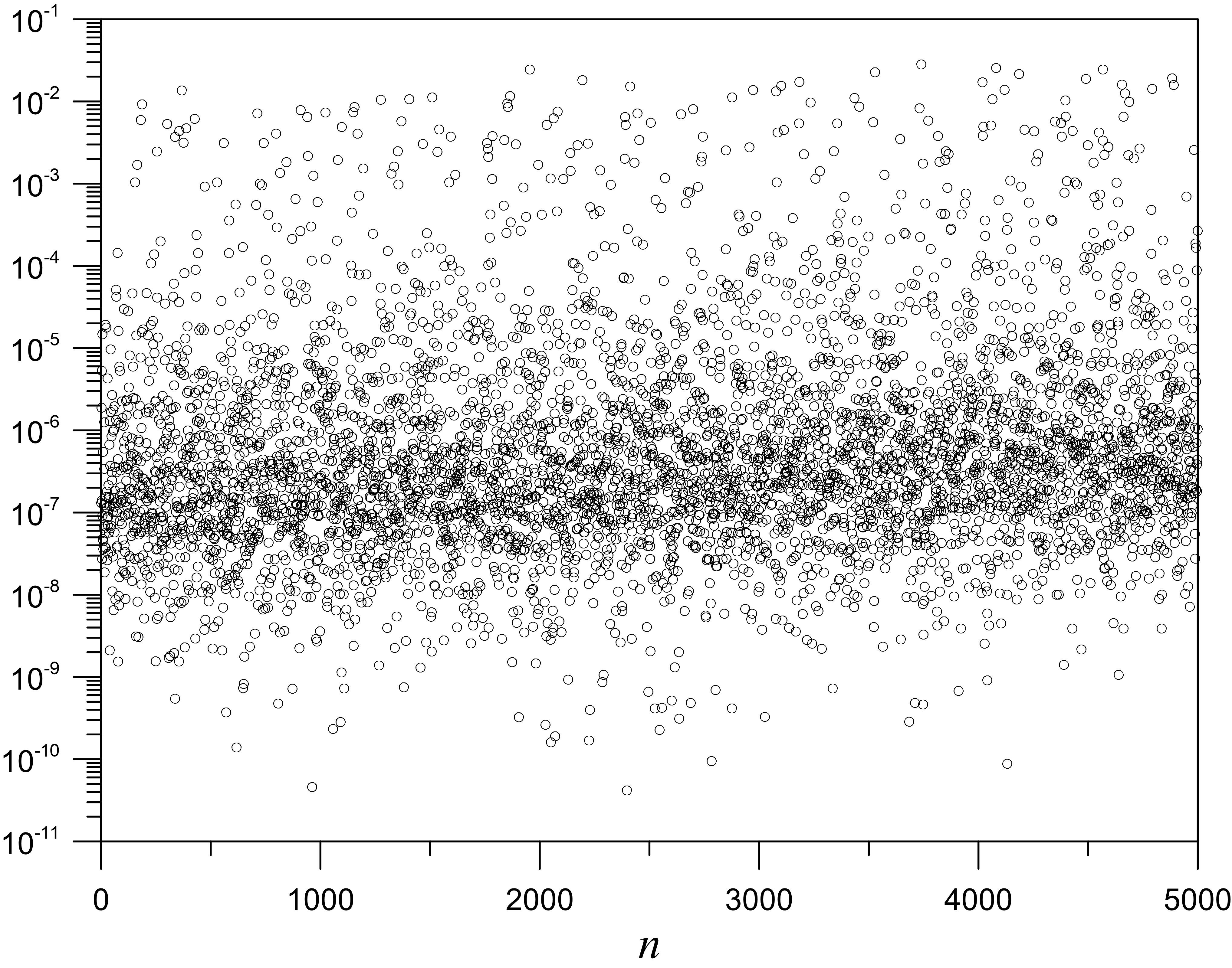} \\
\vspace{-1.5cm}
\vspace{1.7cm}Fig.5 The  plot of the  closest  to the
Khnichin constant $K_0$  values of  the ``running''  geometrical  means  $K_n(m)$.   \\
\end{center}
\end{figure}

Let the rational $P_k/Q_k$ be the $n$-th  partial convergent of the continued fraction:
\bee
\frac{P_k}{Q_k}=[a_0; a_1, a_2, a_3, \ldots, a_k].
\eee
For almost all real numbers $r$ the denominators of the  finite continued fraction
approximations fulfill  \cite[chapter III,  \S 15] {Khinchin}:
\bee
\lim_{k\rightarrow \infty} \big(Q_k(r)\big)^{1/k} = e^{\pi^2/12\ln2} \equiv L_0 = 3.275822918721811\ldots
\eee
where $L_0$ is called the  Khinchin---L{\'e}vy's constant  \cite[\S 1.8]{Finch}.
Again the set of exceptions to  the above limit is  of the Lebesgue measure zero and it includes rational numbers,
quadratic irrational etc.

Let the rational $P_{l(n)}(\gamma_n)/Q_{l(n)}(\gamma_n)$ be the $l$-th
partial convergent of the continued fractions (\ref{g-cfr})  of $\gamma_n$:
\bee
\frac{P_{l(n)}(\gamma_n)}{Q_{l(n)}(\gamma_n)}={\bf a}(n) \doteq\gamma_n .
\eee
For each Stieltjes  constant $\gamma_n$ we   calculated the partial convergents  $P_{l(n)}(\gamma_n)/Q_{l(n)}(\gamma_n)$
using the recurrence:
\[
P_0=a_1, ~~Q_0=1, ~~P_1=1+a_1a_2, ~~Q_1=a_1
\]
\bee
P_k=a_kP_{k-1}+P_{k-2}, ~~Q~_k=a_kQ_{k-1}+Q_{k-2}, ~~k\geq 2.
\eee
Next from these denominators $Q_{l(n)}(\gamma_n)$ we have calculated the quantities
$L_n(l(n))$:
\bee
L_n(l(n))= \left(Q_{l(n)}\right)^{1/l(n)}, ~~~n=1, 2, \ldots , 5000.
\label{def-L}
\eee
The obtained values of $L_n(l(n))$  are presented in the Fig.6.
These  values scatter around the red line representing the Khinchin---L{\'e}vy's
constant $L_0$  and are contained in the interval $(L_0-0.053, L_0+0.053)$.   Again this plot is  somehow misleading
because there are  Stieltjes constant $\gamma(n)$  for which there appear  sign changes of  $L_0-L_n(m), ~m=1,2,\ldots, l(n)$.
As in the case of $K_n(m)$  Fig.7 presents the number of
sign changes of the difference $L_n(m)-L_0$ of the denominator of the $m$-th convergent $P_m/Q_m$
\bee
S_L(n)={\rm number ~ of ~ such ~ {\it m}~ that} ~~~(L_n(m+1)-L_0)(L_n(m)-L_0)<0.
\eee
The maximal number of sign  changes was 922  and appeared for the Stieltjes constant $\gamma_{771}$ and there were 117 gammas
without sign changes.

Finally we looked into the distribution of the values of partial quotients $a_l(n)$.  The  Gauss--Kuzmin theorem \cite[chapter III,  \S 15] {Khinchin}
asserts that the density $d(k)$  of  the denominators $a_m$,  $m=1, 2,  \ldots l$,  with the value $k$ is given by
\bee
\lim_{l \to \infty} \frac{\sharp\{m: a_m=k\}}{l}=\log_2\Big(\frac{1+\frac{1}{k}}{1+\frac{1}{1+k}}\Big)
\label{Gauss-Kuzmin}
\eee
for almost  all  real  numbers.  In the Fig. 9 the results are presented for the Stieltjes  constants.

\begin{figure}[h]
\vspace{-0.3cm}
\begin{center}
\includegraphics[width=0.8\textwidth, angle=0]{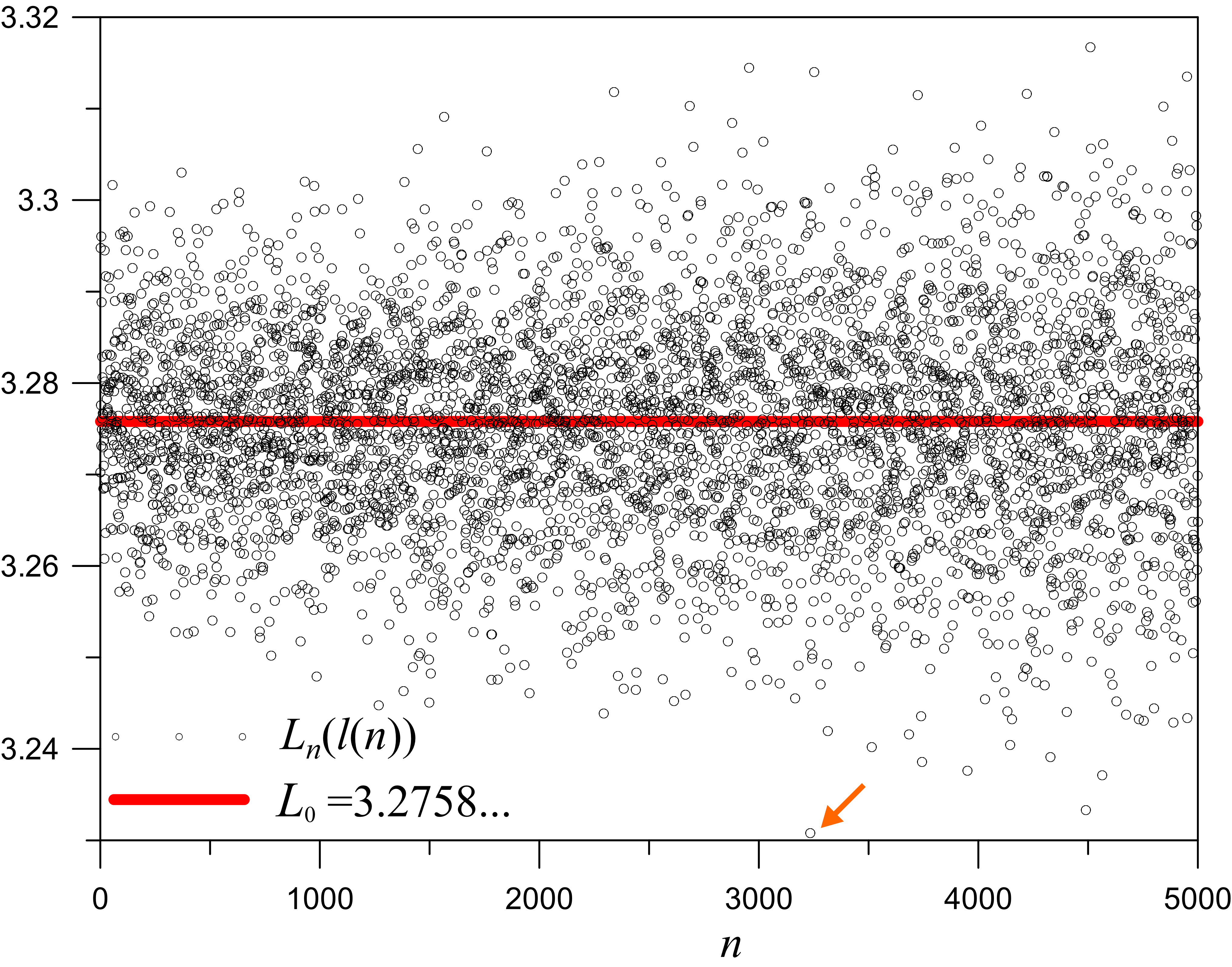} \\
\vspace{0.0cm}
Fig.6  The plot of $L_n(l(n))$ for $n=1,2,3, \ldots, 5000$.  There are 352 points  closer  to $L_0$ than $0.001$ and 38
closer  to $L_0$ than $0.0001$ The largest  value of  $|L_0-L_n(l(n))|$ is $4.503\times 10^{-2}$ and it occurred for the Stieltjes
constant number $l=3235$ (marked with the red  arrow),  the smallest  value of  $|L_0-L_n(l(n))|$ is  $2.336\times 10^{-6}$
and it occurred  for the Stieltjes constant number $l=3226$.   \\
\end{center}
\end{figure}

\begin{figure}[h]
\vspace{-0.3cm}
\begin{center}
\includegraphics[width=0.8\textwidth, angle=0]{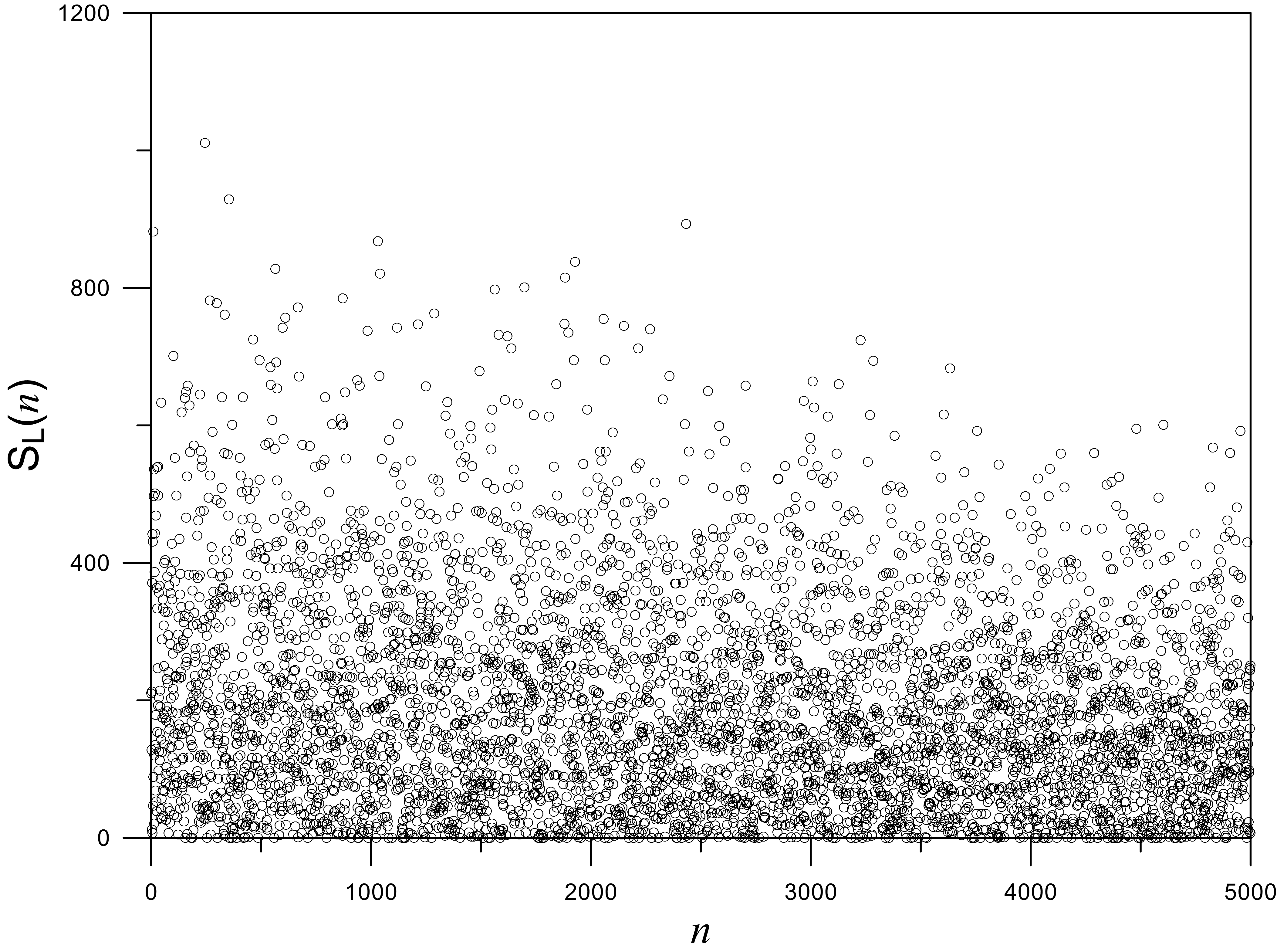} \\
\vspace{0.0cm}
Fig.7  The number of sign changes $S_L(n)$  for each  $n$,  i.e.  the number of  such $m$ that $(L_n(m+1)-L_0)(L_n(m)-L_0)<0$
(the initial transient values  of $m$ were skipped---  sign changes were detected for $m=100, 101, \ldots l(n)$).  \\
\end{center}
\end{figure}




\begin{figure}[h]
\begin{center}
\includegraphics[width=0.8\textwidth, angle=0]{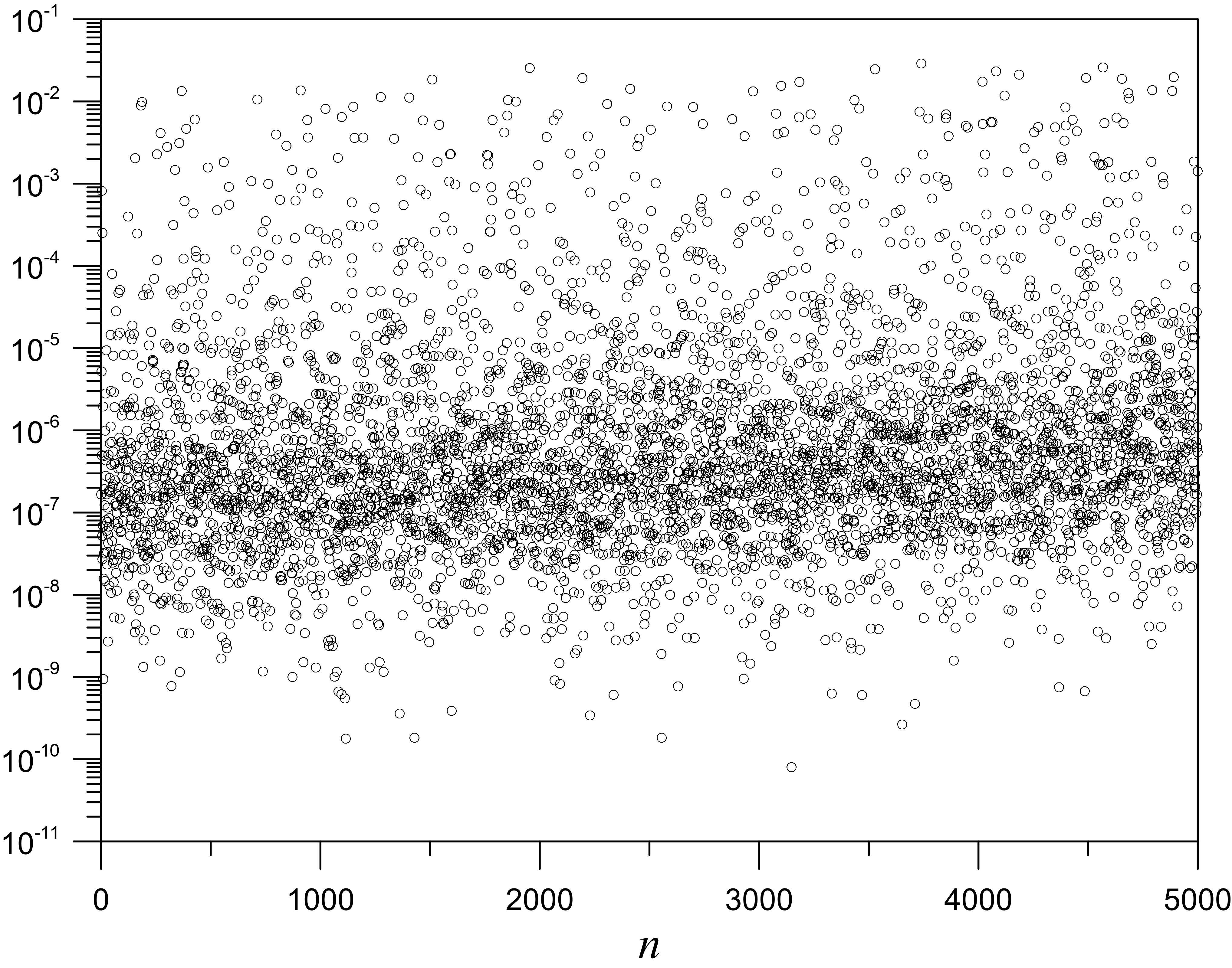} \\
\vspace{-1.5cm}
\vspace{1.7cm}Fig.8 The  plot of the  closest  to the
Khinchin--L{\'e}vy  constant $L_0$  values of  the ``running''  values of  $\sqrt[m]{Q_n(m)}, n=0, 1,2, \ldots, 5000$.   \\
\end{center}
\end{figure}

\begin{figure}[h]
\begin{center}
\includegraphics[width=0.8\textwidth, angle=0]{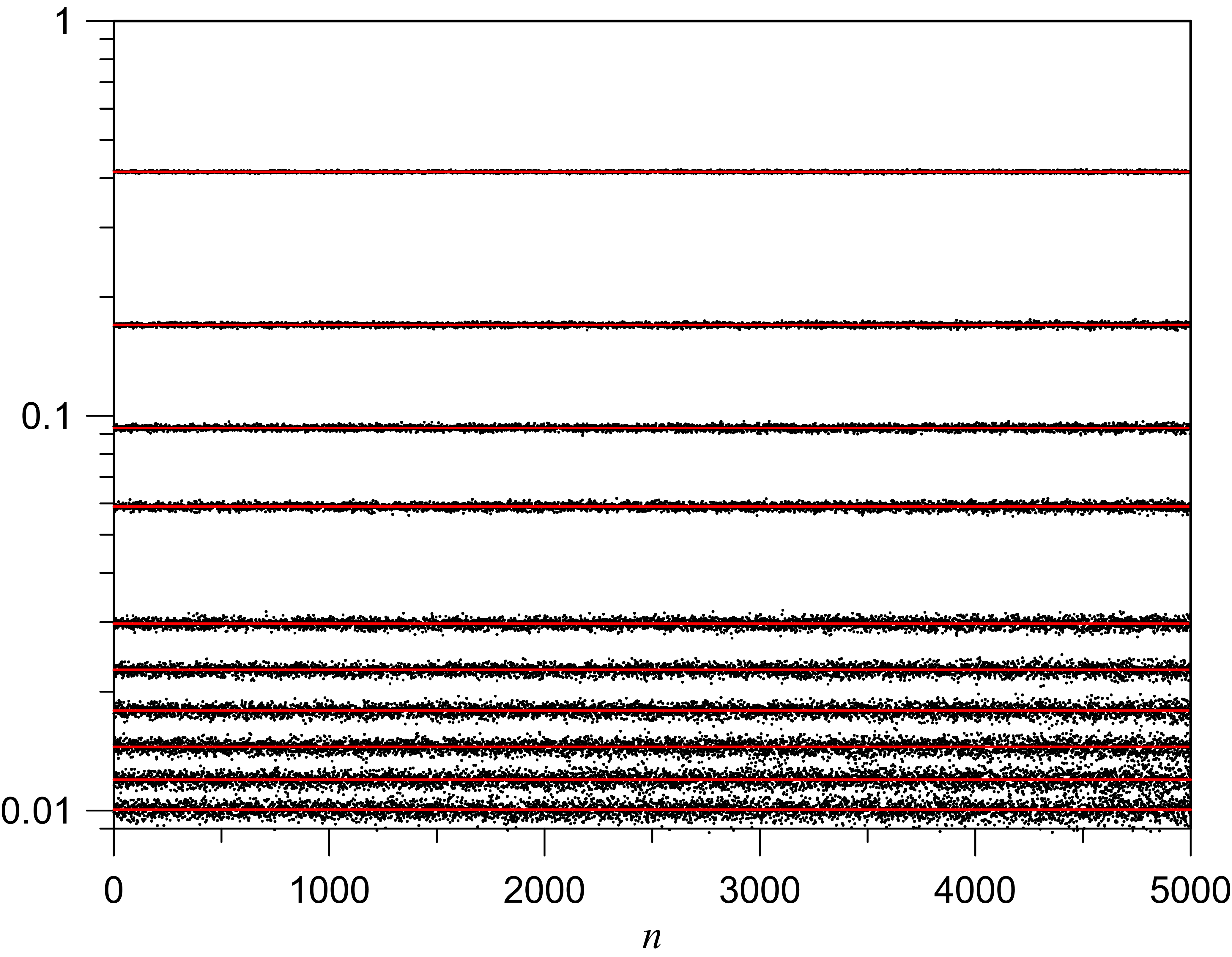} \\
\vspace{-1.5cm}
\vspace{1.7cm}Fig.9 The  plot of the  density of partial quotients $a_l$  equal to $k=1, 2, \ldots, 10 $ from top  to bottom
for  first  5000  Stieltjes  constants. In red   are  the  values of \eqref{Gauss-Kuzmin}  plotted. The $y$  axis is logarithmic to  move the
plots apart. \\
\end{center}
\end{figure}

\section{Appendix: Obtaining high precision numerical values of Stieltjes
constants}

In 1997 it was shown by one of the authors  of the present note \cite{Maslanka}   (M.K.)   that the
Riemann zeta function may be expressed as%
\begin{eqnarray}
\zeta (s) &=&\frac{1}{s-1}\left[ A_{0}+\left( 1-\frac{s}{2}\right)
A_{1}+\left( 1-\frac{s}{2}\right) \left( 2-\frac{s}{2}\right) \frac{A_{2}}{2!%
}+...\right] =  \label{Hyper} \\
&=&\frac{1}{s-1}\sum\limits_{k=0}^{\infty }\frac{A_{k}}{k!}%
\prod\limits_{i=1}^{k}\left( i-\frac{s}{2}\right) = \\
&=&\frac{1}{s-1}\sum\limits_{k=0}^{\infty }\frac{\Gamma \left( k+1-\frac{s}{%
2}\right) }{\Gamma \left( 1-\frac{s}{2}\right) }\frac{A_{k}}{k!}\qquad s\in
\mathbb{C}
\backslash \{1\}
\end{eqnarray}%
where%
\begin{eqnarray}
A_{k} &=&\sum\limits_{j=0}^{k}\left( -1\right) ^{j}\binom{k}{j}(2j+1)\zeta
(2j+2)=  \label{Ak} \\
&=&\frac{1}{2}\sum\limits_{j=0}^{k}\binom{k}{j}(2j+1)\frac{\left( 2\pi
\right) ^{2j+2}B_{2j+2}}{\left( 2j+2\right) !}
\end{eqnarray}
Here  $B_n$  denotes the $n^{{\rm th}}$  Bernoulli  numbers.    
However, the particular choice of nodes in $s=2,4,6,...$, albeit the most
natural, is by no means the only one. One only requires that the prescribed
points be strictly equally spaced. For the purpose of present calculations
we choose the following sequence of points:%
\begin{equation*}
1+\varepsilon ,1+2\varepsilon ,1+3\varepsilon ,...
\end{equation*}%
where $\varepsilon $ is certain real, not necessarily small number.

More precisely, define certain entire function $\varphi $ as:%
\begin{equation*}
\varphi (s):=(s-1)\zeta (s)\qquad s\neq 1
\end{equation*}%
together with $\varphi (1)=1$ which stems from the appropriate limit. Then,
instead of (\ref{Hyper}), we have%
\begin{equation*}
\varphi (s)=\sum\limits_{k=0}^{\infty }\frac{\Gamma \left( k-\frac{s-1}{%
\varepsilon }\right) }{\Gamma \left( -\frac{s-1}{\varepsilon }\right) }\frac{%
\alpha _{k}}{k!}
\end{equation*}%
with%
\begin{equation}
\alpha _{k}=\sum\limits_{j=0}^{k}\left( -1\right) ^{j}\binom{k}{j}\varphi  (1+j\varepsilon )
\label{Ak-general}
\end{equation}
Note that coefficients $\alpha _{k}$\ depend on $\varepsilon $ but we shall
for simplicity drop temporarily this dependence in notation.

As mentioned in the    Introduction   the Stieltjes constants are essentially
coefficients of the Laurent series expansion of the zeta function around its
only simple pole at $s=1$:%
\begin{equation}
\zeta (s)=\frac{1}{s-1}+\sum\limits_{n=0}^{\infty }\frac{\left( -1\right)
^{n}}{n!}\gamma _{n}\left( s-1\right) ^{n}  \label{ZetaExpansion2}
\end{equation}

Now directly from (\ref{ZetaExpansion}) we have:%
\begin{equation*}
\gamma _{n}=\left. \frac{(-1)^{n}}{n+1}\frac{d^{n+1}}{ds^{n+1}}\varphi(s)\right\vert _{s-1}.
\end{equation*}%
Then, after some elementary calculations, we get the following useful result:%
\begin{equation}
\gamma _{n}=\frac{(-1)^{n}n!}{\varepsilon ^{n+1}}\sum%
\limits_{k=n+1}^{\infty }\frac{(-1)^{k}}{k!}\alpha _{k}S(k,n+1)
\label{gamma-gen}
\end{equation}%
where $S(k,i)$ are signed Stirling numbers of the first kind. Note that in
the literature there are different conventions concerning denotation and
indices of Stirling numbers which can be confusing. Here we shall adopt the
following convention involving the Pochhammer symbol:%
\begin{equation*}
\left( x\right) _{k}\equiv \frac{\Gamma (k+x)}{\Gamma (x)}%
=\prod\limits_{i=0}^{k-1}(x+i)=(-1)^{k}\sum%
\limits_{i=0}^{k}(-1)^{i}S(k,i)x^{i}
\end{equation*}

Denoting%
\begin{equation*}
\beta _{nk}\equiv (-1)^{n+k}\frac{n!}{k!}\frac{S(k,n+1)}{\varepsilon ^{n+1}}
\end{equation*}%
we can rewrite (\ref{gamma-gen}) as formally an infinite matrix product%
\begin{equation}
\gamma _{n}=\sum\limits_{k=n+1}^{\infty }\beta _{nk}\;\alpha _{k}
\label{gamma-gen 1}
\end{equation}%
The summation over $k$ starts from $n+1$ since $\beta _{nk}\equiv 0$ for $%
k\leq n$.\ Accuracy of $\alpha _{1}$ is equal to accuracy of precomputed
values of $\varphi(s)$ in equidistant nodes. When $k$ grows the accuracy of
consecutive $\alpha _{k}$ quickly tends do zero. Thus there always exists
certain cut-off value of $k=k_{0}$. Therefore the summation in (\ref{gamma-gen 1}) may be performed to this value:%
\begin{equation}
\gamma _{n}=\sum\limits_{k=n+1}^{k_{0}}\beta _{nk}\;\alpha _{k}
\label{gamma-gen 2}
\end{equation}%
(Numerical experiment confirm that  adding more terms do not affect the  value of the sum \eqref{gamma-gen 2}.)
As pointed earlier $\varepsilon $ need not to be small, however, choosing
smaller $\varepsilon $ greatly accelerates convergence of the series.
However, it also turns out that smaller $\varepsilon $ implies smaller $k_{0}
$. What is really important: All significant digits of $\gamma _{n}$
obtained from the finite sum (\ref{gamma-gen 2}) are correct.

Of course, $\gamma _{n}$ eventually does not depend on $\varepsilon $   
although $\alpha _{k}$ as well as the rate of convergence of (\ref{gamma-gen})
does. In fact series (\ref{gamma-gen}) converges for any value of $%
\varepsilon >0$\ but  the rate of convergence becomes terribly small for $%
\varepsilon \gg 1$. On the other hand, the smaller $\varepsilon $ the faster
the rate of convergence. However, since $\alpha _{k}$ also depends on $%
\varepsilon $, choosing smaller value for $\varepsilon $ requires higher
accuracy of precalculated values of $\varphi(s)$ which in turn may be very time
consuming. Hence, an appropriate compromise in choosing $\varepsilon $ is
needed.

Formula (\ref{gamma-gen}) is particularly suited for numerical calculations.
As already pointed above, one has to choose parameter $\varepsilon $ in
order to optimally perform the calculations. Typically the algorithm has
three simple steps:\bigskip

1. Tabulating   $\varphi(1+j\varepsilon ),j=0,1,2,...$ This requires appropriate
choosing of parameter $\varepsilon $ (see below) and is most time consuming.
The most convenient for this seems small but extremely efficient program
PARI/GP which has implemented particularly optimal zeta procedure.    The first of authors   used
Cyfronet ZEUS computer in Cracow,  where calculating single value of $\varphi(s)$ with $51000$
significant digits requires about 13 minutes.   Since this procedure may easily be
parallelized therefore in order to  compute 10000 values of $\varphi$ 20 independent routines were performed
(each calculating 500 values of $\varphi$) which  took nearly one week.

2. Calculating $\alpha _{k}$ using \eqref{Ak-general} and the precomputed
values.

3. Calculating Stieltjes constants using \eqref{gamma-gen}.
\bigskip

(Contrary to the above step 1 which requires a powerful computer, steps 2
and 3 can be quickly performed on a typical PC.) Several properties
concerning accuracies may be obtained experimentally. It should be stressed
out that given $\alpha _{k}$ calculating single $\gamma _{n}$ with accuracy
of about $50000$ digits requires several minutes on a very modest PC machine.

{\bf  Acknowledgement}: One of the authors (KM)
would like to express his gratitude to the
Academic Computer Center Cyfronet, AGH,
Cracow, for the computational grant of 1000
hours under the PL-Grid project (Polish
Infrastructure for Supporting Computational
Science in the European Research Space).


\end{document}